\documentclass[10pt,a4paper,fleqn]{article}
\usepackage{amsmath}
\usepackage{amssymb}
\usepackage{amscd}
\usepackage[dvips]{graphicx}

\newtheorem{theorem}{Theorem}

\newtheorem{corollary}[theorem]{Corollary}
\newtheorem{definition}[theorem]{Definition}
\newtheorem{example}[theorem]{Example}
\newtheorem{lemma}[theorem]{Lemma}
\newtheorem{proposition}[theorem]{Proposition}
\newtheorem{remark}[theorem]{Remark}
\newtheorem{conjecture}[theorem]{Conjecture}

\newenvironment{proof}{{\bf Proof. }}{\hfill$\rule{1ex}{1ex}$\par\medskip}

\begin{document}
\newcommand{\bt}{\begin{theorem}}
\newcommand{\et}{\end{theorem}}
\newcommand{\bd}{\begin{definition}}
\newcommand{\ed}{\end{definition}}
\newcommand{\bs}{\begin{proposition}}
\newcommand{\es}{\end{proposition}}
\newcommand{\bp}{\begin{proof}}
\newcommand{\ep}{\end{proof}}
\newcommand{\be}{\begin{equation}}
\newcommand{\ee}{\end{equation}}
\newcommand{\ul}{\underline}
\newcommand{\br}{\begin{remark}}
\newcommand{\er}{\end{remark}}
\newcommand{\bex}{\begin{example}}
\newcommand{\eex}{\end{example}}
\newcommand{\bc}{\begin{corollary}}
\newcommand{\ec}{\end{corollary}}
\newcommand{\bl}{\begin{lemma}}
\newcommand{\el}{\end{lemma}}
\newcommand{\bj}{\begin{conjecture}}
\newcommand{\ej}{\end{conjecture}}

\title{Free Steiner triple systems and their automorphism groups}

\author{A. Grishkov, D. Rasskazova, M. Rasskazova, I. Stuhl}

\date{}
\footnotetext{\par}

\maketitle

\begin{abstract}
 The paper is devoted to the study of free objects in the variety
 of Steiner loops and of the combinatorial
 structures behind them, focusing on their
 automorphism groups.

 We prove that all automorphisms
 are tame and the automorphism group is not finitely
 generated if the loop is more than $3$-generated.
 For the free Steiner loop with $3$ generators we describe the generator elements of the automorphism group and some relations between them.
\end{abstract}

\section{Introduction}

Steiner triple systems as special block designs are a major part of combinatorics, and
there are many interesting connections developed between these combinatorial structures
and their algebraic aspects. In this paper we consider Steiner triple systems from
algebraic point of view, i.e., we study the corresponding Steiner loops.
Diassociative loops of exponent $2$ are commutative, and the variety of all
diassociative loops of exponent $2$ is precisely the variety of
all Steiner loops, which are in a one-to-one correspondence with
Steiner triple systems (see \cite{GP}, p. 310).

Since Steiner loops form a variety (moreover a Schreier variety),
we can deal with free objects. Consequently, we use the term {\it free Steiner triple systems} for the combinatorial objects
corresponding to free Steiner loops. A summary of results about varieties
of Steiner loops, Steiner quasigroups and free objects in the varieties can be found in \cite{E}.

We give a construction of free Steiner loops, determine their multiplication groups
(which is a useful knowledge for loops, see \cite{NS}, Section 1.2)). This problem
for finite Steiner loops was partly solved in \cite{SS} and in \cite{SS1} in the case of finite oriented Steiner loops. We also show that the nuclei of the free Steiner loops are trivial,
 which is an indicator of how distant these loops are from groups.

The automorphism group of a Steiner triple system $\mathfrak{S}$ coincides with the
automorphism group of the Steiner quasigroup as well as with the automorphism group of
the Steiner loop associated with $\mathfrak{S}$. Any finite group is the automorphism group of a Steiner triple
system (\cite{M}, Theorem 8, p. 103).
This motivated the goal of our paper to study automorphisms of the free Steiner triple systems.
We prove that ${\rm{(i)}}$ all automorphisms of the free Steiner loops are tame, and ${\rm{(ii)}}$ the automorphism
group of a free Steiner loop is not finitely generated when the loop is generated by more than $3$ elements.

We also determine the generators of the automorphism group of the $3$-generated free Steiner loop and give conjectures about automorphisms of this loop. These conjectures fit the context of the work \cite{U} on linear Nielsen-Schreier varieties of algebras.

\section{Preliminaries}

A set $L$ with a binary operation $L\times L\longrightarrow L:(x,y)\mapsto x\cdot y$ is
called a {\it loop}, if for given $a, b$, the equations $a\cdot y = b$ and
$x\cdot a = b$ are uniquely solvable, and there is an element $e\in L$ such that $e\cdot x = x\cdot e =
x$ for all $x\in L$. A loop is called {\it diassociative} if every two elements generate
a group.

A loop $L$ is called {\it Steiner loop} if $x\cdot (x\cdot y)=y$ holds
for all $x, y\in L$ and $x^{2}=e$ for all $x\in L$, where $e$ is the identity of $L$.

A {\it Steiner triple system} $\mathfrak{S}$ is an incidence
structure consisting of points and blocks such that every two
distinct points are contained in precisely one block, and any block
has precisely three points. It is a well-known fact that a Steiner triple system
of order $m$ exists if and only if $m\equiv 1, 3\, \hbox{mod} \,(6)$
(cf. \cite{Pf}, Definition V.1.9).

To a given Steiner triple system, there correspond two different constructions
leading distinct algebraic structures.

A Steiner triple system $\mathfrak{S}$
determines a multiplication on the pairs of different points $x, y$ taking as a product the
third point of the block joining $x$ and $y$.
Defining $x\cdot x=x$ we get a {\it Steiner quasigroup}
associated with $\mathfrak{S}$. Adjoining an element $e$ with $ex=xe=x$, $xx=e$ we obtain
the Steiner loop $S$.

Conversely, a Steiner loop $S$ determines a Steiner triple system
whose points are the elements of $S\setminus \{e\}$, and
the blocks are the triples $\{x, y, xy\}$ for $x\neq y\in S\setminus \{e\}$. The quasigroup or loop obtained in this way
is called an {\it exterior} Steiner quasigroup
or an {\it exterior} Steiner loop.
This yields the first of the aforementioned constructions. Because it is more popular
than the other one, the term 'exterior' will be omitted.

To describe the second construction, let $a\in \mathfrak{S}$ be some fixed element and $IS=(\mathfrak{S},a,\cdot)$ be a main isotope of
the quasigroup associated to $\mathfrak{S}$ via the multiplication $x\cdot y=y\cdot x=(ax)(ay)$.
Then $x^2=x\cdot x=(ax)(ax)=ax,$ and hence $x^2\cdot y^2=xy,$ $x^3=x(ax)=a$ and
$(xy)y=(x^2\cdot y^2)^2\cdot y^2=x$.

Conversely, from a commutative loop $S$ with identities
$x^3=1, (x^2y^2)^2y^2=x,$ a Steiner triple system can be recovered, with blocks $\{x,y,x^2y^2\}$
and $\{a,x,x^2\}$ for any $x\neq y\neq a$. This construction in a different framework appears in
\cite{E} p 23.
A loop obtained in this way is called an {\it interior Steiner loop}.

The {\it left, right}, respectively,  {\it middle nucleus} of a loop $L$ are the sub\-groups of
$L$ defined by
\[N_l(L) = \{u; \;(u\cdot x)\cdot y = u\cdot (x\cdot y), \; x, y\in L\},\]
\[N_r(L) = \{u; \;(x\cdot y)\cdot u = x\cdot (y\cdot u), \; x, y\in L\},\]
\[N_m(L) = \{u; \;(x\cdot u)\cdot y = x\cdot (u\cdot y), \; x, y\in L\}.\]
The intersection $N(L) = N_l(L)\cap N_r(L)\cap N_m(L)$ is the {\it nucleus} of $L$.

The {\it commutant} $C(L)$ of a loop $L$ is the subset consisting of all elements $c\in L$
such that $c\cdot x = x \cdot c$ for all $x\in L$. The {\it center} $Z(L)$ of $L$ is
the intersection $C(L)\cap N(L)$.

For any $x\in L$ the maps $\lambda_{x}:y\mapsto x\cdot y$ and
$\rho_{x}:y\mapsto y\cdot x$ are the {\it left} and the {\it right
translation}s, respectively. The permutation group generated by the left and right translations
of loop $L$ is called the {\it multiplication group} of $L$, and the stabilizer of
the neutral element is called the {\it inner mapping group} of $L$. These basic facts can be found
in \cite{Ch2}.

\section{Free Steiner loops}
Constructions of free Steiner loops have been given by several authors: see e.g., \cite{E},
\cite{MS}. Nevertheless, we provide here a specific construction; it will help to incorporate a
 transparent interpretation and to establish a natural system of notation.

Let $\tt X$ be a finite ordered set and let $W(\tt X)$ be a set of non-associative $\tt X$-words.
The set $W(\tt X)$ has an order such that
$v>w$ if and only if $|v|>|w|$ or $|v|=|w|>1$, $v=v_1v_2,$ $w=w_1w_2,$ $v_1>w_1$ or
$v_1=w_1,v_2>w_2$. Next, we define
the set $S({\tt X})^{\ast}\subset W({\tt X})$ of $S$-words by induction on the length of word:
\begin{itemize}
  \item ${\tt X}\subset S({\tt X})^{\ast}$,
  \item $wv\in S(\tt X)^{\ast}$ precisely if, $v, w\in S(\tt X)^{\ast}$, $|v|\leq |w|$, $v\neq w$ and if
  $w=w_1\cdot w_2$, then $v\neq w_i$, $(i=1, 2)$.
\end{itemize}
\noindent
On $S({\tt X})=S({\tt X})^{\ast} \cup \{\emptyset\}$ we define a multiplication in the following manner:
\begin{enumerate}
  \item $v\cdot w=w\cdot v=vw$ if $vw\in S({\tt X})$,
  \item $(vw)\cdot w=w\cdot (vw)=w\cdot(wv)=(wv)\cdot w=v$,
  \item $v\cdot v=\emptyset$.
\end{enumerate}
A word $v(x_1, x_2, ... ,x_n)$ is {\it irreducible}, if $v\in S({\tt X})^{\ast}$.

\bs
The set $S({\tt X})$ with the multiplication as above is a free Steiner loop with free generators $\tt X$.
\es

\bp
The definition implies that $S({\tt X})$ is commutative, $a\cdot(a \cdot b)=b$ for all $a, b\in S({\tt X})$
and if $a, b\in S({\tt X})$ then $\{a, b, a\cdot b, \emptyset\}$ is a group of order $4$ and exponent $2$.
Hence $S({\tt X})$ is free diassociative of exponent $2$, i.e., a free Steiner loop.
\ep

Let $(G, H, B)$ be a Baer triple (see \cite{B}), i.e., $G=BH$ is a group, $H$ is a subgroup in $G$, $B$ is a set of representatives for $G/H$ with
$b^2=1, b\in B, B\cap H=1$, and for any $b_1, b_2 \in B$ there exists $b_3\in B$ such that
$b_1b_2=b_3h_1$, $b_2b_1=b_3h_2$, where $h_1, h_2\in H$.

Define a multiplication on $B$ by $b_1\ast b_2 = b_3 = b_2\ast b_1$. Clearly $b\ast b=1$ and
$(b_1\ast b_2)\ast b_2=b_1$. Indeed,
$(b_1\ast b_2)\ast b_2= b_3\ast b_2= b_2\ast b_3 = b_2\ast (b_2b_1{h_2}^{-1})=b_1$ since
$b_2b_2b_1{h_2}^{-1}=b_1{h_2}^{-1}$. This yields that $(B, \ast)$ is a Steiner loop.
We call such a decomposition $G=BH$ an {\it S-decomposition}.

If the intersection $\operatornamewithlimits{\cap}\limits_{x\in G}H^{x}=\{1\}$ then $G\simeq
{\rm{Mult}}((B, \ast))$.

We note that any Steiner loop
can be constructed in the above fashion. Indeed, let $G={\rm{Mult}}(B)$ be the multiplication group of the
Steiner loop $B$ and let $B_0=\{R_b|b\in B\}$, $H=<R_aR_bR_{ab}|a, b\in B>$.
Then $G=B_0H$ is an S-decomposition.

\bs
Let $G={\rm{Mult}}(S({\tt X}))$ be the group of right multiplications of the free Steiner loop $S({\tt X})$. Then
\begin{enumerate}
  \item $G=\operatornamewithlimits{\ast}\limits_{v\in S({\tt X})^{\ast}}C_v$ is a free product of cyclic groups of order $2$;
  \item $G$ acts on $S({\tt X})$, and $G=\{R_v| v\in S({\tt X})\}{\rm{Stab}}_G(\emptyset)$. Moreover, the inner mapping group ${\rm{Stab}}_G(\emptyset)$
  is a free subgroup of $G$ generated by $R_vR_wR_{vw}$, $v, w\in S({\tt X})$.
\end{enumerate}
\es

\bp
The claims follow from the consideration above and from the
definition of the free product of groups
or this fact can be found in \cite{Sm} Sec. 11.3.

The subgroup ${\rm{Stab}}_G(\emptyset)$ is free by
the Kurosh subgroup theorem \cite{K} p. 17.
\ep

\bs
If $x$, $y$ are different elements of the free Steiner loop $S({\tt X})$ and $|{\tt X}|>2$,
then there is an element $z\in S({\tt X})$ such that
$$
(xy)z\neq x(yz).
$$
\es

\bp
Let $x = v_1(x_1,...,x_n)$ and $y = v_2(x_1,...,x_n)$. Suppose we choose the element $z$ in the shape
$z=v_2(x_1,...,x_n)\cdot x_j$, where $x_j$ is one of the generators different from the last letter
of $v_2(x_1,...,x_n)$. Then we have that
$$
(xy)z = (v_1(x_1,...,x_n)\cdot v_2(x_1,...,x_n))(v_2(x_1,...,x_n) x_j)
$$
$$
\neq v_1(x_1,...,x_n) \cdot (v_2(x_1,...,x_n)\cdot v_2(x_1,...,x_n) x_j) = v_1(x_1,...,x_n) x_j =
x(yz).
$$
\ep

As was mentioned earlier, the nucleus of a loop can be interpreted as a 'measure' of the non-associativity. As a
corollary of the previous Proposition, we can conclude that the free Steiner loops are
'very far' from groups:

\bc
The nucleus and therefore the center of free Steiner loops are\\ trivial.
\ec

\section{Automorphisms}

Let ${\tt Y}=\{y_1, y_2,..., y_n\}$ be a set of free generators of $S({\tt X})$. Then $\varphi : {\tt Y}\longrightarrow S({\tt X})$,
$\varphi(y_1)=y_1\cdot v$, $\varphi(y_i)=y_i$, $(i=2,...,n)$, $v\in S({\tt Y}\setminus y_1)$ is
an automorphism of $S({\tt X})$, called an {\it elementary automorphism} (or an {\it ${\tt Y}$-elementary automorphism}) and we will denote it by $\varphi=e_i(v)$.
Let ${\rm{T}}({\tt X})$ denote a subgroup of the group of automorphisms ${\rm{Aut}} (S({\tt X}))$ of loop $S({\tt X})$ generated by
the ${\tt X}$-elementary automorphisms. Automorphisms contained in ${\rm{T}}({\tt X})$ are called ${\it tame}$ (or ${\it {\tt X}-tame}$). In Theorem \ref{tame} below we show that
${\rm{Aut}} (S({\tt X}))={\rm{T}}({\tt X})$.

Let ${\tt Y}=\{y_1, y_2,..., y_m\}\subset S({\tt X})$, then set ${\tt Y}$ is said to be {\it reducible}, if there exist $i$ and
$v\in S({\tt Y}\setminus y_i)$ such that $|y_i\cdot v|<|y_i|$.

Let $S({\tt Z})$ be a free Steiner loop with free generators ${\tt Z}=\{z_1, ..., z_m\}$, let
${\tt Y}=\{y_1,\dots,y_m\}$ be a set of elements of $S({\tt X})$ and let $\varphi:S({\tt Z})\longrightarrow S({\tt Y}):z_i\mapsto y_i$ be
a homomorphism. A set ${\tt Y}$ is called {\it free isometric},
if $\varphi$ is an isomorphism and $|\varphi(v(z_1, ..., z_m))|=||v(z_1, ..., z_m)||$. Here
$||v(z_1, ..., z_m)||$ is the length with weights $|y_1|$, ..., $|y_m|$, it means that $||z_i||=|y_i|$.

\bs
A set ${\tt Y}$ is irreducible if and only if ${\tt Y}$ is free isometric.
\es

\bp
Let ${\tt Y}$ be an irreducible subset of $S({\tt X})$, $S({\tt Z})$ be a free Steiner loop with free generators ${\tt Z}=\{z_1, ..., z_m\}$ and
let $\varphi:S({\tt Z})\longrightarrow S({\tt Y}):z_i\mapsto y_i$ be a homomorphism. We show that $\varphi$ is an
isometric isomorphism.

Let us choose $v\in {\rm{Ker}} \varphi$ of minimal length and set $v=v_1\cdot v_2$, then $\varphi(v_1)=\varphi(v_2)$.
Assume that $v_1=w_1\cdot w_2$ and $v_2=w_3\cdot w_4$ are irreducible, then we have
$\varphi(w_1)\cdot\varphi(w_2)=\varphi(w_3)\cdot\varphi(w_4)$. Suppose that these decompositions are irreducible. Then we get that $\varphi(w_4)=\varphi(w_1)$ or $\varphi(w_4)=\varphi(w_2)$. This yields a contradiction
with the minimality of the choice of $v$ in both cases.

Now we assume, that the decomposition $\varphi(w_1)\cdot\varphi(w_2)$ is reducible, then
$\varphi(w_1)= u_1\cdot\varphi(w_2)$. Hence $u_1= \varphi(w_5)$ and $\varphi(w_1)=\varphi(w_5\cdot w_2)$.
Since the decomposition $w_1\cdot w_2$ is irreducible, $w_1 = w_6\cdot w_7$, $\varphi(w_6)=w_5$ and $\varphi(w_7)=\varphi(w_2)$. Moreover, $w_7\neq w_2$ and therefore $w_2\cdot w_7\in {\rm{Ker}} \varphi$ and $|w_2\cdot w_7|>|v_1\cdot v_2|$. But since $|v_1|>|w_2|$, we have $|v_2|<|w_7|<|w_1|<|v_1|$. This proves the assertion.
\ep

\bc\label{c2}
If ${\tt Y}$ is irreducible then $S({\tt Y})=S({\tt X})$ precisely if ${\tt Y}={\tt X}$.
\ec

Later on we will prove that all automorphisms of the free Steiner loops are tame.

\bt\label{tame}
Let $S({\tt X})$ be a free Steiner loop with free generators ${\tt X}$. Then ${\rm{Aut}} (S({\tt X}))={\rm{T}}({\tt X})$.
\et

\bp
Let $\varphi$ be an automorphism of $S({\tt X})$ and let ${\tt Y}=\varphi({\tt X})$.
We prove that $\varphi\in {\rm{T}}({\tt X})$ by induction on $|{\tt Y}|=\sum_{i=1}^{n}|y_i|$.

First we note that the permutations of ${\tt X}$ are tame automorphisms.
For any transposition $(ij)\in \mathcal{S}_n({\tt X})$ we have $(ij)=\phi\psi\phi$ with
$$\phi=e_i(x_j) \quad \hbox{and} \quad \psi=e_j(x_i).$$
Since the symmetric group $\mathcal{S}_n({\tt X})$ of permutations of $\tt X$ is generated by transpositions, one has
$\mathcal{S}_n({\tt X})\subset {\rm{T}}({\tt X})$.

If $|{\tt Y}|=n$ then $\varphi\in \mathcal{S}_n({\tt X})$ and therefore $\varphi\in {\rm{T}}({\tt X})$.
Now suppose that $|{\tt Y}|>n$. By Corollary \ref{c2} the set ${\tt Y}$ is reducible and hence for
some $i$ and $v=v(y_1, ..., \widehat{y_i}, ..., y_n)$ we
have $|y_i\cdot v|<|y_i|$. By the induction assumption the map
$\psi(x_1,\dots,x_n)=(y_1,\dots,y_{i-1},y_i\cdot v,\dots,y_n)$ induces an ${\tt X}$-tame
automorphism of $S({\tt X}).$ Set
$$w=v(y_1, ..., \widehat{y_i}, ..., y_n)^{\psi^{-1}}=v(x_1, ..., \widehat{x_i}, ..., x_n).$$
Then $\lambda(x_1, ..., x_i, ..., x_n)=(x_1, ..., x_i\cdot w, ..., x_n)$ is an ${\tt X}$-elementary automorphism.
Then $\varphi=\lambda\psi$ since $x_j^{\lambda\psi}=x_j^{\psi}=y_j$ for $j\not=i$ and
$x_i^{\lambda\psi}=(x_i\cdot w)^{\psi}=(y_i\cdot v)\cdot w^{\psi}=(y_i\cdot v)\cdot v=y_i$.

Consequently, $\varphi\in {\rm{T}}({\tt X})$; this completes the proof of the theorem.
\ep

\bl\label{l2}
Let $\phi=e_i(v),$ $v\in S({\tt X}\setminus i)$ be an ${\tt X}$-elementary automorphism and suppose
$u=u_1u_2$ is an ${\tt X}-$irreducible decomposition of a word $u\in S({\tt X})$. Then either $u_1^{\phi}u_2^{\phi}$ is an ${\tt X}-$irreducible decomposition of $u^{\phi}$ or
$u^{\phi}=x_i$, in which case $u_1=x_i$ and $u_2=v$.
\el

\bp
We will use induction in the length $|u|$ of the word $u$.
First suppose that $u_1^{\phi}u_2^{\phi}$ is an ${\tt X}$-reducible decomposition of $u^{\phi}$.
It means that $u_1^{\phi}=u_3u_2^{\phi}$ is also an ${\tt X}$-irreducible decomposition, and hence
$u_1=u_3^{\phi}u_2$. If $u_1=u_3^{\phi}u_2$ is an ${\tt X}$-irreducible decomposition
then $u=u_1u_2$ is ${\tt X}$-reducible, which yields
a contradiction.

Therefore, $u_1=u_3^{\phi}u_2$ is ${\tt X}$-reducible, where
$u_1=x_i$, $u_2=v$, $u_3=x_i$.
Suppose $|u_3|>1$, $|u_2|>1$, $u_3=wx_i$ and $u_2=yv$ $(w\neq x_i\neq v\neq y)$. Then
$u_3^{\varphi}u_2 = [w(x_iv)]\cdot yv$ is ${\tt X}$-reducible if and only if
$w=yv$ or $x_i=y$. In the first case we get that $u_1^{\varphi}u_2^{\varphi}=x_i\cdot yv$ is
${\tt X}$-irreducible. In the second case $u_1=u_3^{\phi}u_2=w$ and
$u_1u_2=w(x_iv)$ is ${\tt X}$-irreducible. Hence, $u_1^{\varphi}u_2^{\varphi}=wx_i$ is
also ${\tt X}$-irreducible decomposition of $u^{\phi}$.
\ep

Define a normal chain of characteristic (${\rm{Aut}} (S({\tt X}))-$invariant) subloops of $S({\tt X})$:
\be\label{eq2}
S_0=S({\tt X})>S_1>S_2>\dots >S_i>\dots\;\;.
\ee
Here $S_0/S_1$ is a group, and for any $i$, \;$Z_i=S_i/S_{i+1}$ is the
center of the factor loop $S_0/S_{i+1}$. Moreover, each $S_i$ is a minimal subloop with these properties.

\vskip 10pt

Now we deal with the question whether the automorphism group of a free Steiner loop with $n$ generators is finitely generated for $n>3$.

\bt
The automorphism group ${\rm{Aut}} (S({\tt X}))$ of the free Steiner loop $S({\tt X})$ is not
finitely generated when $|{\tt X}|>3$.
\et

\bp
Owing to Theorem \ref{tame} and by a discussion afterwards, the group $G={\rm{Aut}} (S({\tt X}))$ is
generated by $\{e_i(v)|v\in S({\tt X})\}$. If $G$ is finitely generated then $G$ is generated
by a set $P=\{e_{j_i}(v_i)|v_i\in S({\tt X}),i=1,\dots,m\}$.

Let $S({\tt X})>S_1>S_2>\dots >S_i>\dots$ be a chain of normal characteristic subloops as in Eqn (\ref{eq2}).
Choose a number $p$ such that $v_i\not\in S_p,$ $i=1,\dots,m$, and $1\not=v\in S_p.$
We assume that
$$e_j(v)=e_{j_1}(v_1)\cdot\dots\cdot e_{j_m}(v_m).$$
For any $w\in S({\tt X})$ we set ${\rm{(i)}}$ $||w||=(s,t),$ if $w=w_1w_2$ is
an ${\tt X}-$irreducible decomposition, $w_1>w_2\not=1,$
$w_1\in S_s\setminus S_{s+1},$ $w_2\in S_t\setminus S_{t+1}$ and ${\rm{(ii)}}$ $||x||=(1,0),$ if $x\in {\tt X}$.

We prove that $||x_je_{j_1}(v_1)\cdot\dots\cdot e_{j_r}(v_r)||=(1,s),$
with $s<m$, by induction in $r$. For $r=1$ this is clear, and we suppose that for $r$ this fact is true.

Set:
$$u=x_1e_{j_1}(v_1)\cdot\dots\cdot e_{j_r}(v_r)e_i(w), \hskip 15 pt
q=x_1e_{j_1}(v_1)\cdot\dots\cdot e_{j_r}(v_r)=q_1q_2.$$
By the induction hypothesis we have $||q_1||=1,$ $||q_2||=s<m.$

If $q_1^{e_i(w)}q_2^{e_i(w)}$ is an ${\tt X}-$irreducible decomposition then
$||q^{e_i(w)}||=||q||=(1,s),$ since $S_s$ is a characteristic subloop. If $q_1^{e_i(w)}q_2^{e_i(w)}$
 is an ${\tt X}-$reducible decomposition then, by Lemma \ref{l2}, $q^{e_i(w)}=u=x_i$ and $||u||=(1,0).$
We obtain that $x_j^{e_j(v)}=x_jv$ and $||x_je_{j_1}(v_1)\cdot\dots\cdot e_{j_m}(v_m)||=(1,s)$, with $s<m$. However, this contradicts to the fact $||x_jv||=(1,m)$. This completes the inductional step.

Therefore our assumption that $G$ is finitely generated does not hold.
\ep
\medskip

The group ${\rm{Aut}} (S({\tt X}))={\rm{T}}({\tt X})$ is generated by ${\tt X}$-elementary automorphisms
$e_i(v)$, $v\in S({\tt X}\setminus\{i\})$, with $e_i(v)^2=1$; this follows from the definition.
Thus, a natural question arises:

\medskip
\noindent
{\bf  Problem 1.}
Which relations exist between ${\tt X}$-elementary automorphisms
of the free Steiner loop $S({\tt X})$?
\medskip

\noindent
We stress, that there is no relation between the elements
$\{e_i(v) | v\in S({\tt X}\setminus x_i) \}$.
\medskip

In what follows we focus on the $3$-generated
free Steiner loop $S(x_1,x_2,x_3)$. Contrary to the case of the automorphism group of free Steiner loop with $n>3$-generators, we prove that the group
${\rm{Aut}} (S(x_1,x_2,x_3))$ is generated
by three involutions $(12),(13)$ and $\varphi=e_1(x_2)$.

\bt\label{t3}
Let $S({\tt X})$ be a free Steiner loop with free generators ${\tt X}=\{x_1, x_2, x_3\}$.
Then the group of automorphisms ${\rm{Aut}} (S({\tt X}))$ is generated by
the symmetric group $\mathcal{S}_3$ and by the elementary automorphism $\varphi=e_1(x_2)$.
\et

\bp
Let $G_0$ be the subgroup of ${\rm{Aut}} (S({\tt X}))$ generated by $\mathcal{S}_3$ and $\varphi$.
If $G_0$ is a proper subgroup, then let $\phi$ be an element of ${\rm{Aut}} (S({\tt X}))\setminus G_0$.
The length of $\phi(x_1, x_2, x_3)=(u, v, w)$ is the sum of the length of the generators under $\phi$,
i.e., $|\phi|=|u|+|v|+|w|$.

The claim of Theorem \ref{t3} can be verified by induction on the length of element $\phi\in {\rm{Aut}} (S({\tt X}))\setminus G_0$.
For $|\phi|=3$, it is trivial. Now if $|\phi|>3$ then by the induction hypothesis we have that if $|\psi|<|\phi|$ then $\psi\in G_0$. By Corollary \ref{c2} the collection $\{u,v,w\}$ is reducible, and we can suppose that
$u=u_0\cdot u_1,$ $u_1\in \{v,w,v\cdot w\}$. There is an automorphism $\alpha$ such that
$\alpha(x_1, x_2, x_3)=(u_0, v, w)$; $\alpha\in G_0$ since $|\alpha|<|\phi|$.
If $u=u_0\cdot v$ then $\phi=\varphi\alpha$. Further, if $u=u_0\cdot(v\cdot w)$ then
\be\label{eq1}
\phi=(13)\varphi (123) \varphi (132) \varphi (13)\alpha.
\ee
Finally, if $u=u_0\cdot w$
then $\phi=(23)\varphi (23)\alpha$. In all three cases $\phi$ is contained in the group $G_0$;
this implies the assertion of the theorem.
\ep

Theorem \ref{t3} implies
\bc
Let $S({\tt X})$ be the Steiner loop with free
generators ${\tt X}=\{a, b, c\}$. Let $Q$ be the stabilizer ${\rm{Stab}}_{{\rm{Aut}} (S({\tt X}))}(c)$
of element $c$ in the automorphism group of $S({\tt X})$. Then
$$Q=<\varphi,\tau,\xi >$$
with
$$
\varphi(a,b,c)=(ab,b,c),\hskip 15 pt \xi(a,b,c)=(ac,b,c),\hskip 15 pt \tau(a,b,c)=(b,a,c).
$$
\ec

\bp
Denote by $Q_0$ the subgroup of $Q$ generated by $\xi, \varphi, \tau$ and
let $\lambda\in Q$ be the map $\lambda(a,b,c)=(v,w,c)$, with $|\lambda|=|v|+|w|$.
Suppose that for every $\gamma\in Q$ with $|\gamma|<|\lambda|$, $\gamma$ is
contained in $Q_0$.

Since $(v,w,c)$ is reducible, we have three possibilities:
$v=v_0w,$ $v=v_0c$ or $v=v_0(wc)$.

Consider the map $\lambda_0(a,b,c)=(v_0,w,c)$; it is contained in
$Q_0$ by induction because $|\lambda_0|< |\lambda|$.

In the first case $\lambda=\varphi\lambda_0$. In the second case
$\lambda=\xi\lambda_0$, and for the mapping $\phi(a,b,c)=(a(bc),b,c)$ we have by Eqn $(\ref{eq1})$ (see the proof of Theorem \ref{t3}).
In the third case $\phi=\tau\xi\varphi\tau\varphi\xi\tau\in Q_0$.

In each case $\lambda\in Q_0$; this fact yields that $Q_0=Q$.
\ep

Let us return to Problem 1. As was mentioned in the proof of Theorem \ref{t3},
any transposition of the symmetric group $\mathcal{S}_n({\tt X})$ on ${\tt X}$ can be written as a product of ${\tt X}$-elementary automorphisms
$$(ij)=e_i(x_j)e_j(x_i)e_i(x_j).$$
Using this description of translations and the equation
$$(i-1,i)(i,i+1)(i-1,i)=(i,i+1)(i-1,i)(i,i+1),$$
we get that
\vskip 10 pt
$e_{i-1}(x_i)e_i(x_{i-1})e_{i-1}(x_i)e_i(x_{i+1})e_{i+1}(x_i)e_i(x_{i+1})e_{i-1}(x_i)
e_i(x_{i-1})e_{i-1}(x_i)=$
\vskip 10 pt
$e_i(x_{i+1})e_{i+1}(x_i)e_i(x_{i+1})e_{i-1}(x_i)e_i(x_{i-1})e_{i-1}(x_i)e_i(x_{i+1})
e_{i+1}(x_i)e_i(x_{i+1}).$
\vskip 10 pt
\noindent
This yields the relation
$$(e_i(x_j)e_j(x_i))^3=1.$$
In the proof of Theorem \ref{t3} we showed a further relation
\vskip 10 pt
$e_1(x_2\cdot x_3)=(13)\varphi (123) \varphi (132) \varphi (13)=$
\vskip 10 pt
$e_1(x_3)e_3(x_1)e_1(x_3)e_2(x_1)e_1(x_2)e_1(x_3)e_3(x_1)e_1(x_3)e_1(x_2)e_1(x_3)e_3(x_1)$
\vskip 10 pt
$\cdot e_1(x_3)e_1(x_2)e_2(x_1)e_1(x_3)e_3(x_1)e_1(x_3).$
\vskip 10 pt
\noindent
These facts suggest the following
\bj\label{cj1}
The group ${\rm{Aut}} (S(x_1,x_2,x_3))$ is generated by three involutions $(12),(13)$
and $\varphi=e_1(x_2)$ with relations
$$(12)(13)(12)=(13)(12)(13), \hskip 20 pt (\varphi(12))^3=(\varphi(13))^4=1.$$
\ej

The analysis of computerised calculations shows that if the Conjecture \ref{cj1} is false then some new relations might exist, between the above involutions, of the type
$$\varphi\sigma_1\varphi\sigma_2\cdots \varphi\sigma_n=1.$$
Here $\sigma_i\in S_3=<(12),(13)>.$ Moreover, $\sigma_i\not=(12)$ or $1;$
if $\sigma_i=(13)$ then $\sigma_{i+1}\not=(13)$. Finally, $n>50$ (for $n\leq 50$
new relations were not found).

In paper \cite{U} it has been proved that the automorphism group of a free algebra of an
arbitrary linear Nielsen--Schreier variety is generated by elementary automorphisms with some specific relations (2)--(4) (\cite{U}, pages 210--211).
If Conjecture \ref{cj1} holds, we will have a similar result for the group
${\rm{Aut}}(S({\tt X}))$ of the free Steiner loop $S({\tt X})$ in the case $|{\tt X}|=3.$

\br
If Conjecture \ref{cj1} is true, group ${\rm{Aut}} (S(x_1,x_2,x_3))$ is
the Coxeter group $<(12),(13),\varphi\;| (\varphi(12))^3 = (\varphi(13))^4 = ((12)(13))^3 = 1 >.$
\er

\bj\label{cj2}
$Q=\{\varphi,\tau,\xi |\xi^2=\varphi^2=\tau^2=(\tau\varphi)^3=1\}.$
\ej

\bt
If Conjecture \ref{cj1} is true then Conjecture \ref{cj2} is also true.
\et

\bp
Suppose that Conjecture \ref{cj1} is true but Conjecture \ref{cj2} is not. Then
there exists a non-trivial word $w=w_1\dots w_n$ formed by the letters $\{\tau,\xi,\varphi\}$ such that
$a^w=a,b^w=b$.
Here the "non-trivial" means that $w$ does not contain the subwords
$\varphi\tau\varphi$ and $\xi\tau\xi\tau$.

Applying induction in $n$, assume that for any non-trivial word $v$
constructed from $\{\tau,\xi,\varphi\}$, of length less then $n$,
the corresponding word in $\{\tau,\pi,\varphi\}$ is non-trivial, where $\pi=(23)$.
Observe that $\xi=\pi\varphi\pi$.
Hence, $w_0=w_1\dots w_{n-1}$ is a non-trivial
word in $\{\tau,\pi,\varphi\}$. We focus on the case where no non-trivial word.
The choice $w_{n-1}=\tau$ implies that $w_n\not=\tau$ and $w_{n-2}\not=\tau$. Furthermore, if $w_n=\xi=\pi\varphi\pi$ then
$w=w_1\dots w_{n-2}\tau\pi\varphi\pi$ is a non-trivial word in $\{\tau,\pi,\varphi\}$.
If $w_n=\varphi$ and $w_{n-2}=\xi$ then $w$ is again a non-trivial word in $\{\tau,\pi,\varphi\}$.
Finally, if $w_n=\varphi$ and $w_{n-2}=\varphi$ then $w$ is not a non-trivial word in $\{\tau,\xi,\varphi\}$,
since $w$ contains the subword $w_{n-2}w_{n-1}w_n=\varphi\tau\varphi$.
\ep

Next, we present as a consequence of the preliminary results,
a connection between the groups of automorphisms of
${\rm{(a)}}$ free Steiner quasi\-groups and their associated ${\rm{(b)}}$ the free exterior Steiner loops and free
interior Steiner loops.

\bt\label{t4}
Let $S({\tt X})$ be a free Steiner quasigroup with free generators ${\tt X}$. Let $ES({\tt X})=S({\tt X})\cup e$ and
$IS({\tt X})$ be its corresponding free exterior and interior Steiner loop, respectively.

Then ${\rm{Aut}}(S({\tt X}))={\rm{Aut}}(ES({\tt X}))$ and ${\rm{Aut}}(IS({\tt X}))\simeq {\rm{Stab}}_{{\rm{Aut}} (ES({\tt X}))}(a)$, where $a\in IS({\tt X})$
is the unit element of loop $IS({\tt X})$.
\et

\vskip 10pt
\noindent
{\bf Acknowledgement}
\vskip 10pt
Izabella Stuhl has been supported by FAPESP
Grant - process number 11/51845-5, and expresses her gratitude to IMS, University of S\~{a}o Paulo,
Brazil, for the warm hospitality.

Alexander Grishkov\\
Institute of Mathematics and Statistics\\
University of S\~{a}o Paulo\\
05508-090 S\~{a}o Paulo, SP, Brazil \vspace{1mm}\\
{\it E-mail}: {\it {}grishkov@ime.usp.br}\\

Diana Rasskazova\\
Novosibirsk State University\\
630090 Novosibirsk, Russia \vspace{1mm}\\
{\it E-mail}: {\it {}ivakirjan@gmail.com}\\

Marina Rasskazova\\
Omsk Service Institute\\
644099 Omsk, Russia \vspace{1mm}\\
{\it E-mail}: {\it {}marinarasskazova1@gmail.com}\\

Izabella Stuhl\\
Institute of Mathematics and Statistics\\
University of S\~{a}o Paulo\\
05508-090 S\~{a}o Paulo, SP, Brazil\\
University of Debrecen\\
H-4010 Debrecen, Hungary \vspace{1mm}\\
{\it E-mail}: {\it {}izabella@ime.usp.br}\\

\end{document}